# Invariance principle for stochastic processes with short memory

## Magda Peligrad[1],[*] and Sergey Utev[2]

*University of Cincinnati and University of Nottingham*

**Abstract:** In this paper we give simple sufficient conditions for linear type processes with short memory that imply the invariance principle. Various examples including projective criterion are considered as applications. In particular, we treat the weak invariance principle for partial sums of linear processes with short memory. We prove that whenever the partial sums of innovations satisfy the $L_p$–invariance principle, then so does the partial sums of its corresponding linear process.

## 1. Motivation: linear processes

Encountered in various applications, the linear time series, moving averages, provide a reach class of examples that are widely studied. For a stationary sequence of innovations $(\xi_i)_{i \in Z}$ and a sequence of constants $(a_j)_{j \in \mathbb{Z}}$, the linear process is defined by

$$(1) \qquad X_k = \sum_{j=-\infty}^{\infty} a_j \xi_{k-j}.$$

Of course, one needs to add some conditions in order for the process to be well defined. In the classical time series analysis $(\xi_i)_{i \in Z}$ is assumed i.i.d. with $\mathbb{E}(\xi_0) = 0$ and $\mathbb{E}(\xi_0^2) < \infty$. Then, $X_k$ is well defined when $\sum_{j=-\infty}^{\infty} a_j^2 < \infty$ and the central limit theorem (CLT) holds for $S_n/stdev(S_n)$ where $S_n = \sum_{j=1}^{n} X_j$ [10]. This CLT is not restricted to i.i.d sequences and the theorem was extended in [19, 21] to martingales and other related structures. However, without additional assumptions on the sequence of constants, the CLT for the general linear process cannot be extended to an invariance principle, not even for independent innovations, as pointed out in [15, 23] and in [16]. A related example is given in Proposition 10 below.

To deal with the problem of the invariance principle, several authors imposed certain regularity conditions on the sequence of constants together with restrictions on the dependence structure of the innovations. Various invariance principles on this line are given for example in [23, 24] and also in [15] among others.

In this paper we shall discuss the case when the coefficients $a_j$ are absolutely summable,

$$(2) \qquad \sum_{i \in \mathbb{Z}} |a_i| < \infty.$$

---

[*]Supported in part by a Charles Phelps Taft Memorial Fund grant and NSA grant, H98230-05-1-0066.
[1]Department of Mathematical Sciences, University of Cincinnati, PO Box 210025, Cincinnati, OH 45221-0025, USA.
[2]School of Mathematical Sciences, University of Nottingham, Nottingham, NG7 2RD, UK.

This case is referred to as short memory, or sometimes as short range dependence. For this situation [11] proved that the central limit theorem (CLT) is preserved under the linear transformation, i.e. the CLT for partial sums of innovations $S_n^{(\xi)} = \sum_{j=1}^n \xi_j$, properly normalized, implies the CLT for $S_n$ under practically the same normalization. This raises the question whether the invariance principle is also preserved under the linear transformation with short memory, that is whenever it holds for the innovations then it also holds for the corresponding linear process.

There are numerous papers dealing with this problem for particular time series with short range dependence and dependent innovations. All these results are asserting that if the innovations have a certain dependence structure, (particular classes of mixing, associated sequences, negative associated sequences, martingales, martingale-like sequences and so on) and satisfy the invariance principle, so does the short range dependent linear process.

To give an example of a result of this type, we mention the well known case when the innovations form a stationary martingale difference sequence [4, 9]. Then, the invariance principle holds,

$$(3) \qquad \frac{S_{[nt]}}{\sqrt{n}} \Longrightarrow \eta A W(t) \quad \text{where} \quad A = \sum_{i \in \mathbb{Z}} a_i \,,\; S_k = \sum_{j=1}^k X_j,$$

where $\eta$ is a random variable measurable with respect to the invariant $\sigma$–field of the stationary sequence $(\xi_i)_{i \in Z}$, denoted by $\mathcal{I}$, $W = \{W(t) \,;\, 0 \leq t \leq 1\}$ is a standard Brownian motion independent of the invariant $\sigma$–field $\mathcal{I}$, $[x]$ denotes the integer part of $x$ and $\Longrightarrow$ denotes weak convergence in $D[0,1]$, the space of cadlag functions on $[0,1]$ endowed with the uniform topology.

In this paper we show that for the short range dependent case the dependence structure is not important and, whenever the innovations satisfy the invariance principle the corresponding short range dependent process is also convergent, provided a certain condition is imposed to $E(\max_{1 \leq j \leq n} |\sum_{k=1}^j \xi_k|)$. We also address the question of the $L_p$– invariance principle ($p \geq 1$) that is

$$(4) \qquad \mathbb{E}(\sup_{0 \leq t \leq 1} |\frac{S_{[nt]}}{\sqrt{n}} - \eta A W(t)|^p) \to 0 \text{ as } n \to \infty$$

and we show that if the innovations satisfy an $L_p$– invariance principle so does the linear process.

For dealing with this problem and other related facts, we develop first a general device to compare the linear combination of processes to one of the initial processes. Another general result will allow to prove the weak convergence for a linear combination of processes by studying only a finite sum. The theory is not restricted to real valued stochastic processes and it can be used for random fields and Hilbert and Banach space valued processes. We then apply the general results to study the asymptotic behavior of linear processes with short memory and to extend the celebrated projective criterion the optimality of our results. We shall also point out that, for a more general linear transformation with short memory, the invariance principle for innovations does not imply in general, the invariance principle for the linear process.

In this paper we shall use the following notations. The space $B(T)$ of bounded functions on $T \subset R^d$ is equipped with the supremum norm $\|x\|_T = \sup_{t \in T} |x(t)|$. The notation $\|X\| = \sqrt{E(X^2)}$ stays for the $L_2$–norm. As already specified in (3), the notation $S_n = X_1 + \cdots + X_n$ will be reserved for the partial sums of the key sequence of interest $(X_i)_i$. Throughout the paper, without changing the distributions we shall redefine the innovations on a probability space rich enough to support a standard Brownian motion.

## 2. The general device

A well-known powerful tool in the proofs of weak convergence results consists of approximating the underlying sequence $(W_n)$ by a double indexed sequence $(W_n^{(m)})$, such that for any $m$, the sequence $W_n^{(m)} \to W^{(m)}$ and then, we deduce the existence of the limit of the original sequence $(W_n)$ by studying the limit of $(W^{(m)})$ as $m \to \infty$. (see for example, Theorem 3.2. in [1]).

In this section we shall establish two variants of this type of approximation, aimed to study the linear combination of stochastic processes.

First, we present a coupling type approximation lemma. We consider a process that can be expressed as a linear combination with short memory of stochastic processes, and approximate it by an individual summand. The context is general enough to be applied to processes indexed by countable sets allowing also the treatment of random fields and, just a modification of language will lead us to higher dimensional spaces and operators.

**Lemma 1.** *Let $\psi^{(n)} = \psi^{(n)}(t)$, $t \in T$, $n = 1, 2, \dots$ be a sequence of stochastic processes that admits the following representation*

$$\psi^{(n)} = \sum_{j \in I} a_j U_j^{(n)} \quad \text{with} \quad \sum_{j \in I} |a_j| < \infty, \tag{5}$$

*where $I$ is a countable set, $(a_j)_{j \in I}$ is a real valued sequence and $U_j^{(n)} = U_j^{(n)}(t)$, $t \in T$, $j \in I$, $n = 1, 2, \dots$ is a double sequence of stochastic processes, satisfying the following conditions.*

$$\sup_n \sup_j \mathbb{E}(\|U_j^{(n)}\|_T) < \infty \tag{6}$$

*and, for each pair $i, j \in I$,*

$$\|U_i^{(n)} - U_j^{(n)}\|_T \to 0 \quad \text{in probability as } n \to \infty. \tag{7}$$

*Fix $e \in I$. Then, with the notation $A = \sum_{j \in I} a_j$, we have*

$$\|\psi^{(n)} - A U_e^{(n)}\|_T \to 0 \quad \text{in probability as } n \to \infty.$$

*Moreover, if for a certain $p \geq 1$ we have $\sup_n \sup_j \mathbb{E}(\|U_j^{(n)}\|_T^p) < \infty$ and the convergence in (7) is in $L_p$ then*

$$\mathbb{E}(\|\psi^{(n)} - A U_e^{(n)}\|_T^p) \to 0 \text{ as } n \to \infty.$$

*Proof of Lemma 1.* First we find a sequence $(I_m)_{m \geq 0}$ of subsets of $I$ such that $I_0 = \emptyset$, $I = \cup I_m$, the set $I_m$ contains exactly $m$ elements and $I_m \subset I_{m+1}$. Observe that for any positive integer $m$ and $e \in I$,

$$\psi^{(n)} - A U_e^{(n)} = \sum_{j \in I} a_j (U_j^{(n)} - U_e^{(n)})$$
$$= \Big( \sum_{j \in I - I_m} a_j (U_j^{(n)} - U_e^{(n)}) \Big) + \Big( \sum_{j \in I_m} a_j (U_j^{(n)} - U_e^{(n)}) \Big).$$

Hence, by the triangle inequality,

$$\|\psi^{(n)} - AU_e^{(n)}\|_T \leq \Big( \sum_{j \in I - I_m} |a_j|(\|U_j^{(n)}\|_T + \|U_e^{(n)}\|_T) \Big)$$
$$+ \Big( \sum_{j \in I} |a_j| \Big) \Big( m \max_{j \in I_m} \|U_j^{(n)} - U_e^{(n)}\|_T \Big)$$

and so, for $\epsilon > 0$,

$$P(\|\psi^{(n)} - AU_e^{(n)}\|_T \geq \epsilon) \leq \frac{4}{\epsilon} \Big( \sum_{j \in I - I_m} |a_j| \Big) \sup_n \sup_j \mathbb{E}(\|U_j^{(n)}\|_T)$$
$$+ P\Big( m \sum_{j \in I} |a_j| \max_{j \in I_m} \|U_j^{(n)} - U_e^{(n)}\|_T > \frac{\epsilon}{2} \Big).$$

We let first $n \to \infty$, and notice that by condition (6) the second term vanishes. Then, we let $m \to \infty$ and the first term vanishes because

$$\sum_{j \in I - I_m} |a_j| \to 0 \quad \text{as } m \to \infty.$$

If the convergence in (7) is in $L_p$ and $\sup_n \sup_j \mathbb{E}(\|U_j^{(n)}\|_T^p) < \infty$ then, obviously, the convergence in the conclusion of the last part of this lemma holds also in $L_p$. □

To comment on the conditions used in this lemma, we shall see later on that condition (6) is in a particular context also a necessary condition. To verify condition (7) one needs a certain linear structure of the process $U_j^{(n)}$. In our next lemma, we avoid this condition. Moreover we show that for stochastic processes satisfying condition (5) the asymptotic behavior is then identified by the limiting behavior of roughly a sum of $m$ sequences, where $m$ is arbitrary fixed positive integer.

**Lemma 2.** *Let $\psi^{(n)} = \psi^{(n)}(t)$, $t \in T = [0,1]^d$, $n = 1, 2, \ldots$ be a sequence of stochastic processes satisfying conditions (5) and (6). Let $I = \cup_m I_m$ for some sequence of increasing subsets $I_m$ of $I$ and assume that for each $m$,*

(8)
$$\sum_{j \in I_m} a_j U_j^{(n)} \Longrightarrow Z_m \text{ as } n \to \infty$$

*in $C[0,1]^d$ (or in $D[0,1]^d$ endowed with the uniform norm) where $Z_m$ is a continuous stochastic process. Then, there exists a limiting stochastic process $Z_m \Longrightarrow Z$ as $m \to \infty$ in $C[0,1]^d$ and $\psi^{(n)} \Longrightarrow Z$ as $n \to \infty$ in $C[0,1]^d$.*

*Proof of Lemma 2.* Let

$$\psi_m^{(n)} = \sum_{j \in I_m} a_j U_j^{(n)}$$

Notice that for $h \geq m$,

$$\mathbb{E}(\|\psi_h^{(n)} - \psi_m^{(n)}\|_T) \leq 2 \Big( \sum_{j \in I - I_m} |a_j| \Big) (\sup_n \sup_j \mathbb{E}(\|U_j^{(n)}\|_T)) \to 0$$

as $h \geq m \to \infty$.

Whence, by the Fatou lemma, the sequence of stochastic processes $(Z_m)_{m \in \mathbb{N}}$ satisfies the Cauchy criterion

$$\mathbb{E}(\|Z_m - Z_h\|_T) \to 0 \quad \text{as } m, h \to \infty.$$

Therefore, $(Z_m)_{m\in\mathbb{N}}$ has a limiting process $Z$ that has sample paths in $C[0,1]^d$.

We also have

$$\mathbb{E}(\|\psi^{(n)} - \psi_m^{(n)}\|_T) \leq 2\Big(\sum_{j\in I-I_m}|a_j|\Big)(\sup_n \sup_j \mathbb{E}(\|U_j^{(n)}\|_T)) \to 0 \text{ as } m \to \infty$$

and thus, by Theorem 3.2. in Billingsley (1999), $\psi^{(n)} \Longrightarrow Z$ as $n \to \infty$ which proves the lemma. □

## 3. Examples

### An i.i.d. example.

We begin by showing that somewhat surprisingly, condition (6) is not just a technical condition and it is necessary in some situations.

**Proposition 3.** *Let $(U_j^{(n)})_{j,n\in\mathbb{N}}$ be a double array of i.i.d. non-negative random variables with finite mean, $\mathbb{E}(U_1^{(n)}) < \infty$. The following conditions are equivalent:*

(i) $U_1^{(n)} \to^P 0$ *as* $n \to \infty$ *and* $\sup_n \mathbb{E}(U_1^{(n)}) < \infty$.
(ii) *For any non-negative sequence* $(a_j)_{j\in\mathbb{N}}$

$$\text{if } \sum_{j\in\mathbb{N}} a_j < \infty \quad \text{then} \quad \sum_{j=1}^{\infty} a_j U_j^{(n)} \to^P 0 \text{ as } n \to \infty.$$

*Proof.* Implication (i) → (ii) follows from Lemma 1. Now, to prove implication (ii) → (i), we notice that for the fixed $n$, $(U_j^{(n)})_{j\in\mathbb{N}}$ is an i.i.d. sequence and, the Kolmogorov three series theorem then implies that

$$\text{for any non-negative sequence } (a_j)_{j\in\mathbb{N}} \quad \text{with } \sum_{j\in\mathbb{N}} a_j < \infty$$

we have,

$$\sum_{j=1}^{\infty} \mathbb{E}g(a_j U_j^{(n)}) \to 0 \quad \text{as } n \to \infty.$$

where $g(x) = xI_{(x\leq 1)} + I_{(I>1)}$. Since for any positive $a$

$$\mathbb{E}g(a\xi) = a\int_0^{1/a} P(\xi \geq t)dt$$

then,

$$\sum_{j=1}^{\infty} \mathbb{E}g(a_j U_j^{(n)}) = \sum_{j=1}^{\infty} \int_0^{1/a_j} P(U_1^{(n)} \geq x)a_j dx$$

$$= \sum_{j=1}^{\infty} \int_0^{\infty} P(U_1^{(n)} \geq x)I_{(x\leq 1/a_j)}a_j dx$$

$$= \int_0^{\infty} P(U_1^{(n)} \geq x)t_a(x)dx$$

where
$$t_a(x) = \sum_{j=1}^{\infty} I_{(x \leq 1/a_j)} a_j.$$

On the other hand, for any non-negative non-increasing function $t(x)$ such that $t(x) \downarrow 0$ we can find $a_j \downarrow 0$ such that

$$t(x) \leq t_a(x) \, , \text{ for all } \, x \geq x_0.$$

and so, for any non-negative non-increasing function $t(x)$ such that $t(x) \downarrow 0$ we have

(9) $$\int_0^{\infty} P(U_1^{(n)} \geq x) t(x) dx \to 0 \quad \text{as } n \to \infty$$

Now, we proceed by contradiction and, without loss of generality, assume that

(10) $$\mathbb{E}(U_1^{(n)}) \to \infty \text{ as } n \to \infty.$$

Then, we can find two increasing sequences $n_k \nearrow \infty$, $M_k \nearrow \infty$ and a positive sequence $Q_k \to \infty$ as $k \to \infty$ such that

$$\int_{M_{k-1}}^{M_k} P(U_1^{(n_k)} \geq x) dx \geq Q_k.$$

Now, we take $t(x) = 1/Q_k$ on the half open interval $x \in (M_{k-1}, M_k]$ and $t(x) = 0$ for $x \leq x_0$. Then,

$$\int_0^{\infty} P(U_1^{(n_k)} \geq x) t(x) dx \geq \frac{1}{Q_k} \int_{M_{k-1}}^{M_k} P(U_1^{(n_k)} \geq x) dx \geq 1 \not\to 0 \text{ as } k \to \infty$$

and so (10) contradicts (9). □

### Linear stationary processes

Let $(\xi_i)_{i \in \mathbb{Z}}$ be a stationary sequence with $\mathbb{E}(|\xi_0|) < \infty$ and $\mathbb{E}(\xi_0) = 0$ and let $\mathcal{I}$ be its invariant $\sigma$–field. Define

(11) $$X_k = \sum_{j=-\infty}^{\infty} a_j \xi_{k-j} \quad \text{and assume} \quad \sum_{i=-\infty}^{\infty} |a_i| < \infty.$$

In addition to the notations defined in (3), we remind the notation

$$S_k^{(\xi)} = \sum_{j=1}^{k} \xi_j \, .$$

Our next proposition allows to compare the partial sums of the innovations to the partial sums of the linear process with short memory.

**Proposition 4.** *Assume the representation* (11) *is satisfied and in addition, there is a constant $C > 0$ and a sequence of positive reals $b_n \to \infty$ such that for all $n$,*

(12) $$\mathbb{E}\Big(\max_{1 \leq j \leq n} |S_j^{(\xi)}|\Big) \leq C b_n$$

and

$$b_n^{-1} \max_{1 \leq j \leq n} |\xi_j| \to^P 0 \quad as \ n \to \infty \ . \tag{13}$$

Then,

$$\frac{1}{b_n}(\max_{1 \leq j \leq n} |S_j - AS_j^{(\xi)}|) \to^P 0 \quad as \ n \to \infty \ . \tag{14}$$

If the innovations are assumed in $L_p$, $p \geq 1$, (12) is replaced by $\mathbb{E}(\max_{1 \leq j \leq n} |S_j^{(\xi)}|^p) \leq Cb_n^p$ and the convergence in (13) holds in $L_p$, then

$$\frac{1}{b_n^p} E(\max_{1 \leq j \leq n} |S_j - AS_j^{(\xi)}|^p) \to 0 \quad as \ n \to \infty \ .$$

*Proof.* This proposition follows from Lemma 1 applied to the representation

$$\frac{1}{b_n} S_{[nt]} = \sum_{i=-\infty}^{\infty} a_i U_i^{(n)} \quad \text{where} \quad U_i^{(n)} = \frac{1}{b_n} \sum_{k=1}^{[nt]} \xi_{k-i}$$

so that conditions (5) and (6) follow from (11) and (12), while the convergence in (7) follows from (13). □

From this proposition, we derive

**Theorem 5.** *Assume that representation (11) and condition (12) are satisfied. Moreover assume that the innovations satisfy the invariance principle $b_n^{-1} S_{[nt]}^{(\xi)} \Longrightarrow \eta W(t)$ as $n \to \infty$ where $\eta$ is $\mathcal{I}$–measurable and $W$ is a standard Brownian motion $[0,1]$ independent on $\mathcal{I}$. Then the linear process also satisfies the invariance principle, i.e. $b_n^{-1} S_{[nt]} \Longrightarrow \eta AW(t)$ as $n \to \infty$ .*

*Proof.* Notice that the convergence in probability in (13) follows from the invariance principle $b_n^{-1} S_{[nt]}^{(\xi)} \Longrightarrow \eta W(t)$ as $n \to \infty$ , since the modulus of continuity is convergent to 0. All the conditions in Proposition 4 are then satisfied which imply the conclusion of the theorem. □

From Theorem 5 we easily derive the following useful consequence.

**Corollary 6** ($L_p$-invariance principle). *Assume the representation (11) holds and $p \geq 1$. Then,*

$$If \quad \mathbb{E}(\sup_{0 \leq t \leq 1} |b_n^{-1} S_{[nt]}^{(\xi)} - \eta W(t)|^p) \to 0 \quad then \quad \mathbb{E}(\sup_{0 \leq t \leq 1} |b_n^{-1} S_{[nt]} - \eta AW(t)|^p) \to 0$$

*as $n \to \infty$.*

### Discussion

*Applications*

Theorem 5 and its Corollary 6 work for many dependent structures such as surveyed in [3, 6, 17], Merlevède, Peligrad and Utev (2006). Various invariance principles can be extended from the original sequence to the linear process with short memory.

Here we mention some traditional and also recently developed dependence conditions for innovations whose partial sums satisfy both a maximal inequality and the invariance principle and therefore Theorem 5 applies.

Let us assume that $(\xi_i)_{i \in \mathbb{Z}}$ is a stationary sequence with $\mathbb{E}(\xi_0^2) < \infty$ and $\mathbb{E}(\xi_0) = 0$ and let $\mathcal{F}_a^b$ be a $\sigma$–field generated by $\xi_i$ with indexes $a \leq i \leq b$. For all the structures bellow the family $(max_{1 \leq k \leq n} S_k^2/n)_{n \geq 1}$ is uniformly integrable and the conclusion of Corollary 6 holds with $b_n = \sqrt{n}$ and with $p = 2$.

(i) Hannan [9], see also its extension to Hilbert space in [4]:

$$\sum_{n=1}^{\infty} \|P_0(\xi_n)\| < \infty \text{ and } \mathbb{E}(\xi_0|\mathcal{F}_{-\infty}) = 0 \text{ a.s.}$$

where $P_k(X) = \mathbb{E}(X|\mathcal{F}_k) - \mathbb{E}(X|\mathcal{F}_{k-1})$ is the projection operator.

(ii) Newman and Wright [18]: $(\xi_i)_{i \in \mathbb{Z}}$ is a negatively associated sequence or, positively associated sequence with

$$\sum_{k=1}^{\infty} cov(\xi_k, \xi_0) < \infty.$$

(iii) Doukhan, Massart and Rio [7]:

$$\sum_{k=1}^{\infty} \int_0^{\tilde{\alpha}(k)} Q^2(u) du < \infty.$$

where $Q$ denotes the cadlag inverse of the function $t \to P(|\xi_0| > t)$ and $\tilde{\alpha}(k) = \alpha(\mathcal{F}_0, \mathcal{F}_n^n) = \sup\{|P(A \cap B) - P(A)P(B)| \; ; \; A \in \mathcal{F}_0, B \in \mathcal{F}_n^n\}$ is the strongly mixing coefficient.

(iv) Dedecker and Rio [5]:

$$\mathbb{E}(X_0 S_n|\mathcal{F}_0) \text{ converges in } \mathbb{L}_1.$$

(v) Peligrad and Utev [20], by developing Maxwell and Woodroofe[12]:

$$\sum_{n=1}^{\infty} \frac{\|\mathbb{E}(S_n|\mathcal{F}_0)\|_2}{n^{3/2}} < \infty.$$

(vi) Peligrad, Utev and Wu [22], which guarantees the $L_p$-invariance principle: For $p \geq 2$.

$$\sum_{n=1}^{\infty} \frac{\|\mathbb{E}(S_n|\mathcal{F}_0)\|_p}{n^{3/2}} < \infty.$$

*Remarks*

(a) If $\mathbb{E}(\xi_0^2) < \infty$ and $b_n \geq \sqrt{n}$, then condition (13) automatically holds.
(b) If the sequence of innovations is ergodic then there is a nonnegative constant $\sigma$ such that $\eta = \sigma$ a.s.
(c) The set of indexes $Z$ can be replaced by $Z^d$ where $d$ is a positive integer allowing for the treatment of random fields.
(d) A natural extension is to consider innovations with values in functional spaces that also facilitate the study of estimation and forecasting problems for several classes of continuous time processes [2]. The linear processes are still defined by the formula (1) with the difference that now, the innovations $(\xi_k)_{k \in \mathbb{Z}}$ are Hilbert space $H$-valued random elements and the sequence of constants is replaced by the sequence of bounded linear operators $\{a_k\}_{k \in \mathbb{Z}}$ from $H$ to $H$. In [13], was treated the problem of the central limit theorem for this case under the summability condition

$$\sum_{j=-\infty}^{\infty} \|a_j\|_{L(H)} < \infty,$$

where $\|a_j\|_{L(H)}$ denotes the usual operators. It was discovered that, if this condition is not satisfied, then the central limit theorem fails even for the case of independent innovations. The approach developed in this paper shows that the central limit theorem results stated can be strengthened to the invariance principle (some results in this direction for strongly mixing sequences are established in [14]).

### *Projective criteria*

In this section we apply the general devise in Lemma 2 to derive a non-stationary projective criteria. Let $X$ be an integrable random variable and $\{\mathcal{F}_j, j \in \mathbb{Z}\}$ be a non-decreasing filtration of $\sigma$-fields, that is $\mathcal{F}_j \subseteq \mathcal{F}_i$ for all $j \leq i$. As before, define the projection operator by

$$P_k(X) = \mathbb{E}(X|\mathcal{F}_k) - \mathbb{E}(X|\mathcal{F}_{k-1}).$$

The next proposition gives a linear representation of the process $n^{-1/2} S_{[nt]}$ in terms of a linear combination of processes involving the sequence of projections.

**Proposition 7.** *Let $(X_i)_{i \in \mathbb{Z}}$ be a square integrable centered sequence which is adapted to the non-decreasing filtration $(\mathcal{F}_i)_{i \in \mathbb{Z}}$. Let $\mathcal{F}_{-\infty} = \bigcap_{i \in \mathbb{Z}} \mathcal{F}_i$ Assume that for all $k = 1, 2, \ldots,$*

(15) $$\mathbb{E}(X_k|\mathcal{F}_{-\infty}) = 0 \ a.s.$$

*and for all $k = 1, 2 \ldots$ and $j = 0, 1, 2, \ldots$*

(16) $$\|P_{k-j}(X_k)\| \leq p_j \quad \text{where } p_j > 0 \text{ and } \sum_{j=0}^{\infty} p_j < \infty .$$

*Then, the process $n^{-1/2} S_{[nt]}$ satisfies the representation (5) with*

$$a_i = p_i \text{ and } U_k^{(n)} = n^{-1/2} \sum_{i=1}^{[nt]} P_{k-i}(X_k) p_i^{-1}$$

*and so $(n^{-1/2} S_{[nt]})$ satisfies the invariance principle when (8) holds.*

*Proof.* From conditions (15) and (16) it follows the following martingale difference decomposition

(17) $$X_k = \sum_{i=0}^{\infty} P_{k-i}(X_k)$$

which proves (5).

On the other hand, to check condition (6) we apply the Doob $L_2$–maximal inequality

$$\mathbb{E}\Big(\max_{1\leq j\leq n}\Big|\sum_{k=1}^{j}P_{k-i}(X_k)\Big|^2\Big) \leq 4\sum_{k=1}^{n}\|P_{k-i}(X_k)\|^2 \leq 4(\sqrt{n})^2(p_i)^2$$

which implies condition (6) and completes the proof of this proposition by Lemma 2. □

If we impose more restrictive degrees of stationarity the conclusion can be strengthen.

**Proposition 8.** *Let $(X_i)_{i\in\mathbb{Z}}$ be a centered, uniformly square integrable sequence of random variables, adapted to the non-decreasing filtration $(\mathcal{F}_i)_{i\in\mathbb{Z}}$. Assume that conditions (15) and (16) are satisfied and in addition for each $t \in [0,1]$, and all $m = 0, 1, 2, \ldots$*

$$(18) \qquad \frac{1}{n}\sum_{j=1}^{[nt]}(P_j(S_{j+m-1} - S_{j-1}))^2 \to^P \eta_m B(t) \text{ as } n \to \infty$$

*where $B(t)$ is a non-random non-decreasing function. Then, there is a random variable $\eta$ such that $n^{-1/2}S_{[nt]} \Longrightarrow \eta B(t)W(t)$ where $W$ is a standard Brownian motion independent of $\eta$.*

*Proof.* We employ the martingale decomposition and then, we use a standard result for martingales. By Proposition 7, we know that the limiting behavior of $n^{-1/2}S_{[nt]}$ is determined by the limiting behavior of the partial sum process

$$n^{-1/2}\sum_{k=1}^{[nt]}Y_k \quad \text{where} \quad Y_k = X_k - \mathbb{E}(X_k|\mathcal{F}_{k-m}),$$

where $m$ is a fixed arbitrary positive integer. Notice that $\mathbb{E}(Y_j|\mathcal{F}_{j-m}) = 0$ almost surely and so the following variables are properly defined

$$\theta_k = \sum_{j=k}^{\infty}\mathbb{E}(Y_j|\mathcal{F}_k) = \sum_{j=k}^{k+m-1}\mathbb{E}(Y_j|\mathcal{F}_k) = Y_k + Q_k \quad a.s.$$

In particular, it is easy to see that

$$\mathbb{E}(\theta_k|\mathcal{F}_{k-1}) = -Y_{k-1} + \sum_{j=k-1}^{\infty}\mathbb{E}(Y_j|\mathcal{F}_{k-1}) = -Y_{k-1} + \theta_{k-1}$$

and we derive the following coboundary decomposition

$$S_n^{(Y)} = \sum_{i=1}^{n}Y_i = M_n + (Q_0 - Q_n) \quad \text{where} \quad M_n = \sum_{k=1}^{n}(\theta_k - \mathbb{E}(\theta_k|\mathcal{F}_{k-1})).$$

By construction and the conditions of the proposition,

$$n^{-1/2}\max_{t\in[0,1]}|Q_0 - Q_{[nt]}| \to^P 0 \text{ as } n \to \infty$$

and so, the limiting behavior of the partial sum process $n^{-1/2}S^{(Y)}_{[nt]}$ is determined by the limiting behavior of the normalized discrete time martingale $n^{-1/2}M_{[nt]}$, $t \in [0,1]$ with uniformly square integrable martingale differences

$$\theta_k - \mathbb{E}(\theta_k|\mathcal{F}_{k-1}) = P_k(S_{k+m-1} - S_{k-1}) .$$

Then the proposition follows easily by standard results for the functional central limit theorem for martingales (see [8], or [1]). □

We now define a stationary filtration as in [12], that is we assume that $X_i = g(Y_j, j \leq i)$ where $(Y_i)_{i \in \mathbb{Z}}$ is an underlying stationary sequence. Denote by $\mathcal{I}$ its invariant sigma field and by $(\mathcal{G}_i)_{i \in \mathbb{Z}}$ an increasing filtration of sigma fields $\mathcal{G}_i = \sigma(Y_j, j \leq i)$. For the case when for every $i$, $\xi_i = Y_i$, and $g(Y_j, j \leq i) = Y_i$, then $\mathcal{G}_i$ is simply the sigma algebra generated by $\xi_j$, $j \leq i$. From the above Proposition 8, by using the stationary ergodic theorem we easily derive the stationary projective criterion contained in the next theorem.

We shall derive a class of invariance principles for linear type statistics. The central limit theorem was treated in [19, 21].

**Theorem 9.** *Let $(X_i)_{i \in \mathbb{Z}}$ be a stationary sequence with $\mathbb{E}(X_0) = 0$ and $\mathbb{E}(X_0^2) < \infty$ and stationary filtration ( $\mathcal{G}_i)_{i \in \mathbb{Z}}$. Let $\mathcal{G}_{-\infty} = \bigcap_{i \in \mathbb{Z}} \mathcal{G}_i$. Assume that*

$$\mathbb{E}(X_0|\mathcal{G}_{-\infty}) = 0 \quad a.s. \quad and \quad \sum_{i \geq 1} \|P_0(X_i)\|_2 < \infty .$$

*Then, there exists an $\mathcal{I}$–measurable positive random variable $\eta$ such that for any Lipschitz function $g$,*

$$n^{-1/2} \sum_{i=1}^{[nt]} g(i/n) X_i \Longrightarrow \eta B_g(t) W(t) \text{ as } n \to \infty$$

*where $W$ is a standard Brownian motion independent of $\eta$ and $B_g(t) = \sqrt{\int_0^t g^2(x)dx}$.*

*Proof.* Since $g(x)$ is Lipschitz, therefore bounded, by Proposition 8, we have only to check condition (18) that reduces to establishing the convergence

(19) $$\frac{1}{n} \sum_{j=1}^{[nt]} g^2(j/n)(P_j(S_{j+m-1} - S_{j-1}))^2 \to^P \eta_m \int_0^t g^2(x)dx$$

as $n \to \infty$. Let

$$G(j/n) = g^2(j/n) , \quad \psi_j = P_j(S_{j+m-1} - S_{j-1}).$$

Notice that by the Birkhoff-Khintchine ergodic theorem there exists the limit

$$\frac{1}{n} \sum_{j=1}^{n} \psi_j \to \eta_m \quad a.s., \text{ as } n \to \infty.$$

On the other hand, because the function $G(x) = g^2(x)$ is bounded and continuous, therefore Riemann integrable, we derive

$$\frac{1}{n} \sum_{j=1}^{[nt]} G(j/n) \to \int_0^t g(x)^2 dx \text{ as } n \to \infty,$$

and, in order to prove (19) it is enough to show that

$$(20) \qquad \frac{1}{n}\sum_{j=1}^{[nt]} G(j/n)(\psi_j - \eta_m) \to 0 \quad a.s., \text{ as } n \to \infty.$$

Let

$$U_j = \sum_{i=1}^{j}(\psi_i - \eta_m)$$

and notice that

$$(21) \qquad \sup_j(|U_j|/j) < \infty \quad \text{and} \quad \sup_{j \geq N}(|U_j|/j) \to 0 \text{ as } N \to \infty$$

almost surely.

We easily get the representation

$$\sum_{j=1}^{[nt]} G(j/n)(\psi_j - \eta_m) = G([nt]/n)U_{[nt]} + \sum_{j=1}^{[nt]-1}(G(j/n) - G([j+1]/n))U_j.$$

Clearly for any fixed $t$, $G([nt]/n)U_{[nt]}/n \to 0$ almost surely. On the other hand, since $g$ is Lipschitz, so is $G = g^2$. Therefore there is a constant $K$ such that $|G(j/n) - G([j-1]/n)| \leq K/n$ and thus

$$\left| \frac{1}{n}\sum_{j=1}^{[nt]-1}(G(j/n) - G([j+1]/n))U_j \right| \leq \frac{1}{n}\sum_{j=1}^{n} j|G(j/n) - G([j+1]/n)||U_j/j|$$

$$\leq \frac{1}{n}\sup_j |U_j/j| \left(\sum_{j=1}^{[\sqrt{n}]} \frac{Kj}{n}\right) + \frac{1}{n}\sup_{j \geq \sqrt{n}}|U_j/j| \left(\sum_{j=1}^{n} \frac{Kj}{n}\right)$$

$$\to 0 \text{ as } n \to \infty$$

almost surely, which proves (20) and thus completes the proof of the theorem. □

We notice that for stationary sequences, when the filtration is not stationary the result is not true if we only assume the summability of projective norms. We can have the decomposition in (5) but condition (6) will not be satisfied and the invariance principle will fail.

**Proposition 10.** *There exists a stationary Gaussian positively associated linear process*

$$X_k = \sum_{i=0}^{\infty} t_i Y_{k-i}$$

*where $t_i \geq 0$ and $Y_i$ is an i.i.d. sequence of standard normal Gaussian variables and a non-decreasing filtration $\mathcal{F}_k$ such that*

(i) $X_k$ *is adapted to* $\mathcal{F}_k$.
(ii) $E(X_k|\mathcal{F}_{-\infty}) = 0$ *a.s. for all $k$.*
(iii)
$$\sup_{k \geq 0} \sum_{j \in \mathbb{Z}} \|P_j(X_k)\| < \infty.$$

(iv) $\sigma_n^2 = \mathrm{Var}(S_n)$ *is not regularly varying with index* 1, *more exactly, there exists a positive c and a subsequence* $k_n \to \infty$ *such that* $\mathrm{Var}(S_{k_n}) \geq ck_n^2/\ln^4(k_n)$ .

(v) $\sigma_n^{-1}S_{[nt]}$ *does not satisfy the invariance principle, i.e does not converge to the standard Brownian motion.*

*Proof.* Let $(Y_k)_{k \in \mathbb{Z}}$ be an i.i.d. sequence of standard normal variables. Let $\mathcal{G}_i = \sigma(Y_j, j \leq i)$, $\mathcal{G}_\infty = \sigma(Y_j, j \in \mathbb{Z})$. For a positive integer $r$, let

$$n_r = 4^r \quad \text{and} \quad u_r = (n_{r+1} - n_r)^{-1/2} r^{-4} = 1/(3r^4 2^r).$$

We take

$$t_j = u_r \quad \text{when} \quad n_r < j \leq n_{r+1}, r = 1, 2, \ldots$$

We also take $t_i = 0$ when $i < n_1$. Next, as a filtration we take

$$\mathcal{F}_i = \mathcal{G}_{-n_r} \quad \text{when} \quad -n_{r+1} < i \leq -n_r, r = 1, 2, \ldots$$

and in addition

$$\mathcal{F}_i = \mathcal{G}_\infty \quad \text{when} \quad i \geq -n_1.$$

We notice first that $X_k$ is $\mathcal{F}_k$–measurable and so the sequence $(X_k)_{k \in \mathbb{Z}}$ is adapted to the filtration $(\mathcal{F}_k)_{k \in \mathbb{Z}}$, which proves (i).

Clearly,

$$\mathcal{F}_{-\infty} = \bigcap_{i \in \mathbb{Z}} \mathcal{F}_i = \bigcap_{r \in \mathbb{N}} \mathcal{G}_{-n_r} = \mathcal{G}_{-\infty},$$

therefore,

$$\mathbb{E}(X_k|\mathcal{F}_{-\infty}) = \mathbb{E}(X_k|\mathcal{G}_{-\infty}) = 0 \quad a.s.$$

which proves (ii).

Now, let us compute the projection operator. For $i \geq -n_1$, we have for all $X_k$, since they are $\mathcal{G}_\infty$–measurable,

$$\mathbb{E}(X_k|\mathcal{F}_i) = \mathbb{E}(X_k|\mathcal{G}_\infty) = X_k \quad a.s.$$

implying that, for $i \geq -n_1$,

$$P_i(X_k) = \mathbb{E}(X_k|\mathcal{F}_i) - \mathbb{E}(X_k|\mathcal{F}_{i-1}) = 0 \quad a.s.$$

Now, for $i$ such that $-n_{r+1} < i - 1 < i \leq -n_r$,

$$P_i(X_k) = \mathbb{E}(X_k|\mathcal{F}_i) - \mathbb{E}(X_k|\mathcal{F}_{i-1})$$
$$= \mathbb{E}(X_k|\mathcal{G}_{-n_r}) - \mathbb{E}(X_k|\mathcal{G}_{-n_r}) = 0 \quad a.s.$$

Finally, for $i$ such that $-n_{r+1} = i - 1 < i \leq -n_r$,

$$P_i(X_k) = \mathbb{E}(X_k|\mathcal{F}_i) - \mathbb{E}(X_k|\mathcal{F}_{i-1}) = \mathbb{E}(X_k|\mathcal{G}_{-n_r}) - \mathbb{E}(X_k|\mathcal{G}_{-n_{r+1}})$$
$$= \mathbb{E}\left(\sum_{j=0}^\infty t_j Y_{k-j}|\mathcal{G}_{-n_r}\right) - \mathbb{E}\left(\sum_{j=0}^\infty t_j Y_{k-j}|\mathcal{G}_{-n_{r+1}}\right)$$
$$= \left(\sum_{j=k+n_r}^\infty t_j Y_{k-j}\right) - \left(\sum_{j=k+n_{r+1}}^\infty t_j Y_{k-j}\right)$$
$$= \sum_{j=k+n_r+1}^{k+n_{r+1}} t_j Y_{k-j} \quad a.s.$$

and hence, for such an $i$,

$$\|P_i(X_k)\|_2^2 = \sum_{j=k+n_r+1}^{k+n_{r+1}} t_j^2$$

By construction, the sequence $t_j$ is non-increasing and so

$$\|P_i(X_k)\|_2^2 \leq \sum_{j=n_r+1}^{n_{r+1}} t_j^2 = (n_{r+1} - n_r)u_r^2 = r^{-8}$$

which implies that

$$\sup_{k \geq 0} \sum_j \|P_j(X_k)\|_2 \leq \sum_{r=1}^{\infty} r^{-4} < \infty$$

proving (iii).

To compute the variance, we observe that

$$\mathrm{Var}(S_n) = \sum_{i=-\infty}^{\infty} \Big(\sum_{k=1}^{n} t_{k-i}\Big)^2.$$

So, for $n = n_{r+1}$, $-[n_{r+1}/3] \leq i \leq 0$ and $[n_{r+1}/2] \leq k \leq n_{r+1}$, we have $n_r < k - i$ and therefore

$$\mathrm{Var}(S_{n_{r+1}}) \geq \sum_{i=-[n_{r+1}/3]}^{0} \Big(\sum_{k=[n_{r+1}/2]}^{n_{r+1}} t_{k-i}\Big)^2$$

$$\geq \sum_{i=-[n_{r+1}/3]}^{0} \Big(\sum_{k=[n_{r+1}/2]}^{n_{r+1}} u_r\Big)^2$$

$$\geq (u_r n_{r+1})^2 n_{r+1}/12 = \Big(\frac{4^{r+1}}{3r^4 2^r}\Big)^2 \frac{4^{r+1}}{12} = \frac{4^{2(r+1)}}{9r^8} = \frac{n_{r+1}^2}{9r^8}$$

and so

$$\mathrm{Var}(S_{n_{r+1}}) \geq n_{r+1}^2/(9\log_4(n_{r+1}))$$

which proves (iv). We conclude that the variance is not regularly varying with index 1. It is well known however that the weak convergence of $\sigma_n^{-1} S_{[nt]}$ to $W(t)$, standard Brownian motion, implies that variance is regularly varying with index 1 and then the invariance principle cannot hold. □